\documentclass[a4paper,article,10pt]{memoir}
\usepackage{microtype}
\usepackage{CCLAuthor}
\usepackage[pdftex, pdfauthor={Philippe H. Trinh}]{hyperref}
\usepackage[usenames,dvipsnames]{xcolor}
\definecolor{Maroon}{cmyk}{0, 0.87, 0.68, 0.32}
\hypersetup{colorlinks=true, citecolor=MidnightBlue, linkcolor=MidnightBlue}

\usepackage[reqno]{amsmath}
\usepackage{amssymb}
\usepackage{amsthm}  
\theoremstyle{plain}

\theoremstyle{definition}

\numberwithin{theorem}{chapter}
\usepackage[round,authoryear]{natbib} 
\usepackage[]{graphicx}               
\usepackage{psfrag}
\begin{document}
%
%
%
\title{Exponential Asymptotics and Stokes Line Smoothing for Generalized
Solitary Waves}
%
%
\author{%
    Philippe H. Trinh
    \\ \smallskip\small
    Oxford Centre for Industrial and Applied Mathematics, \\ Mathematical Institute, University of Oxford, UK
    }
    \maketitle
%
%
%
\begin{abstract}
In a companion paper, Grimshaw (\emph{Asymptotic Methods in Fluid Mechanics}, 2010, pp. 71--120) has demonstrated how techniques of Borel summation can be used to elucidate the exponentially small terms that lie hidden \emph{beyond all orders} of a divergent asymptotic expansion. Here, we provide an alternative derivation of the generalized solitary waves of the fifth-order Korteweg-de Vries equation. We will first optimally truncate the asymptotic series, and then smooth the Stokes line. Our method provides an \emph{explicit} view of the switching-on mechanism, and thus increased understanding of the Stokes Phenomenon.
\end{abstract}

\CCLsection{Introduction}

The \emph{Stokes Phenomenon} describes the puzzling event in which exponentially
small terms can suddenly appear or disappear when an asymptotic expansion is
analytically continued across key lines \emph{(Stokes lines)} in the Argand
plane---\emph{``as it were into a mist,''} Stokes once remarked in 
\citeyear{stokes_1902}. 

Fortunately, much of the inherent vagueness of this phenomenon, as well as its
deep implications for the study of asymptotic approximations has been examined
since Stokes' time (see \cite{boyd_1999} for a comprehensive review). In another
paper of this volume by \cite{grimshaw_2011}---henceforth referred to as [Grimshaw]---it was
shown how Borel summation can be used to reveal the exponentially small waves
found in the fifth-order Korteweg-de Vries equation (5KdV).

In this review paper, we will show how the methodology outlined in
\cite{olde_1995} and \cite{chapman_1998} can be used as an alternative treatment
of the 5KdV equation. The procedure is as follows: (1) Expand the solution as a
typical asymptotic expansion, (2) find the behaviour of the late-order terms
($n\to\infty$), and (3) optimally truncate the expansion and examine the
remainder as the Stokes lines are crossed. 

\begin{figure}[htbp]
\centering
\includegraphics[width=1.0\textwidth]{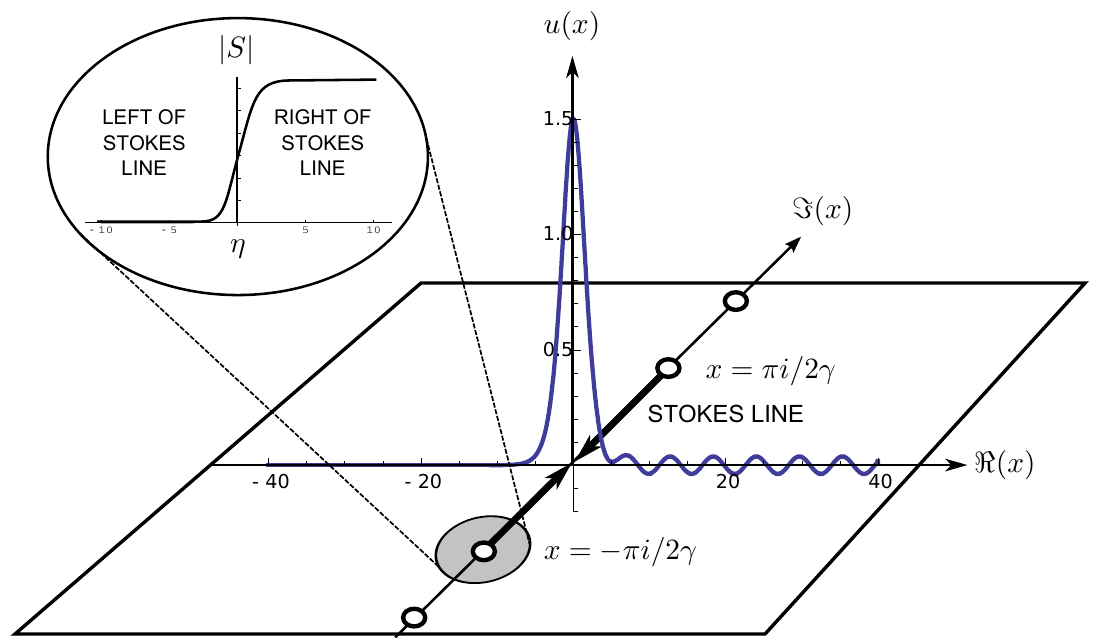}
\caption{The analytic continuation of the traditional asymptotic solution (the
classical solitary wave) contains singularities up and down the imaginary axis,
with Stokes lines emanating from each of these singularities. By re-scaling near
the singularities and optimally truncating, we will be able to observe the
smooth switching-on of the exponentially small terms (top-left). \label{fig}}
\end{figure}

The location of the Stokes lines, as well as the details of the Stokes
Phenomenon and resultant exponentials are intrinsically linked to the late-
order terms of the asymptotic approximation---thus, as we proceed through Steps
1 to 3, we are effectively deriving the beyond-all-orders contributions by
\emph{decoding} the divergent tails of the expansion. The novelty in
this approach (in contrast to the one shown in [Grimshaw]) is that all the
analysis is done in the (complexified) physical space, rather than in
Borel-transformed space. This provides us with a special vantage point---to see
the \emph{smooth} switching-on of the exponentially small terms as each Stokes
line is crossed (see Figure \ref{fig}). Come, let us stare into Stokes' mist.

\vspace*{-1.5\baselineskip}
\CCLsection{Generalized Solitary Waves and the 5KdV} 

We will consider the existence of solutions to the 5KdV equation,
\begin{gather} \label{eq:5kdv}
\epsilon^2 u_{xxxx} + u_{xx} + 3u^2 - cu = 0
\end{gather}
with $u \to 0$ as $x \to \pm \infty$. Although the problem is for
$x \in \mathbb{R}$, it will be important to consider the effects of allowing 
$u$ and $x$ to be complex.

\CCLsubsection{Initial Asymptotic Analysis and Late Terms}

\noindent We begin as usual by substituting the asymptotic expansions,
\begin{equation*}
u = \sum_{n=0}^\infty \epsilon^{2n} u_n \text{\quad and \quad} c =
\sum_{n=0}^\infty \epsilon^{2n} c_n
\end{equation*}

\noindent into Equation (\ref{eq:5kdv}). This yields the first two orders as,
\begin{align}
u_0 &= 2\gamma^2 \textrm{sech}^2(\gamma x) & c_0 &= 4\gamma^2 \label{eq:u0} \\
u_1 &= -10\gamma^2 u_0 + \left(\frac{15}{2}\right)u_0^2 & c_1 &= c_0^2
\end{align}
while at $O(\epsilon^{2n})$,
\begin{equation} \label{eq:Oepn}
u_{(n-1)xxxx} + u_{nxx} + 6u_0u_n - c_0u_n + \ldots = 0.
\end{equation}

\noindent Here, the key observation is that there exists \emph{singularities} in
the analytic continuation of the leading order solution, $u_0(x)$ at $x = \pm
\pi i/2\gamma, \pm 3\pi i/2\gamma, \ldots$ This use of \emph{ill-defined}
approximations in order to represent perfectly well-defined phenomena is one of
the caveats of singular asymptotics, but one would feverishly hope that a
singularity far from the region of interest (in this case, $x \in \mathbb{R}$)
has little effect on the approximations!

Unfortunately this is not the case. We see from Equation (\ref{eq:Oepn}) that at
each order, $u_n$ is partly determined by differentiating $u_{n-1}$ twice and
thus each additional order adds to the power of the singularities in
the early terms. We would therefore expect the late terms of the asymptotic
expansion to exhibit \emph{factorial over power} divergence of the form,
\begin{equation} \label{eq:ansatz}
u_n \sim \frac{Q(z)\Gamma(2n + \gamma)}{[\chi(z)]^{2n + \gamma}}, \quad
\text{as
$n \to \infty$.}
\end{equation}

\noindent Here, $\gamma$ is a constant, while $Q(z)$ and $\chi(z)$ are functions
to be determined. Substituting this ansatz into Equation (\ref{eq:Oepn}) yields
at leading order,
\begin{equation}
-\left( \frac{d\chi}{dz}\right)^4 + \left( \frac{d\chi}{dz}\right)^2 = 0, \quad
\text{as
$n \to \infty$.}
\end{equation}
Now from the above discussion, we would expect that $\chi = 0$ at
the relevant singularities, $x = \sigma_i$ for some $i$; we then conclude that
$\chi' = \pm 1$ and thus without loss of generality, $\chi = x - \sigma_i$. In
general, $u_n$ will be a sum of terms of the form (\ref{eq:ansatz}), one
for each singularity. However along the real axis, the behaviour of $u_n$ will
be dominated by those singularities closest to the axis and thus we need only
concern ourselves with the singularities at $x = \pm \sigma \equiv \pm
i\pi/2\gamma$. Finally, at next order as $n \to \infty$, we find that $Q(z) =
\Lambda$, a constant. 

The determination of $\gamma$, $\Lambda$, and in fact, the Stokes line
smoothing in the next section will require an analysis near each of the two
singularities, $x = \pm \sigma$; for brevity, we will henceforth focus on the
singularity at $x = \sigma$ in the upper-half plane.

First, since by Equation (\ref{eq:u0}), $u_0 \sim -2/(x-\sigma)^2$ as $x \to
\sigma$, we must require that $\gamma = 2$. Second, in order to determine
the final constant $\Lambda$, we need to re-scale near the singularity $x = 
\sigma$, express the leading-order inner solution as a power series (in inner
coordinates) and match with the outer solutions. In the end, however, $\Lambda$
is determined by the numerical solution of a canonical \emph{inner} problem. As
was shown in [Grimshaw], $\Lambda \approx -19.97$. 

Finally, let us discuss the significance of $\chi$. Following Dingle
\citeyearpar{dingle_book}, we expect there to be a Stokes line wherever $u_n$
and $u_{n+1}$ have the same phase as $n \to \infty$, or in this case where,
\begin{equation}
 \Im\left[ -\chi^2 \right] = 0 \text{ and } \Re\left[ -\chi^2 \right] \geq 0.
\end{equation}
Thus there exist Stokes lines from $x = \pi i/2\gamma$
\emph{down} the imaginary axis and from $x = -\pi i/2\gamma$ \emph{up} the
imaginary axis (as illustrated in Figure \ref{fig}). In the next section, we
will optimally truncate the asymptotic expansion and examine the switching-on of
exponentially small terms as these two Stokes lines are crossed.

\CCLsubsection{Optimal Truncation and Stokes Smoothing}

\noindent By now we have entirely determined the late terms of the asymptotic
expansion. In order to identify the exponentially small waves, we truncate the
expansion and study its remainder,
\[ u = \sum_{n=0}^{N-1} \epsilon^{2n} u_n + R_N(x). \]
Substitution into Equation (\ref{eq:5kdv}) yields the equation
\begin{equation} \label{eq:remeq}
\epsilon^2 R_N'''' + R_N'' + 6u_0R_N - c_0R_N + \ldots \sim \epsilon^{2N}u_N'', 
\end{equation}
which, using Stirling's formula, we can write the right-hand side as
\begin{multline} \label{eq:rhs}
 \epsilon^{2N}u_N'' \sim \epsilon^{2N} \frac{\Lambda(-1)^N \Gamma(2N + \gamma +
2) \chi'}{\chi^{2N+\gamma+2}} \\ \sim \epsilon^{2N}
\frac{\Lambda(-1)^N \left\{ \sqrt{2\pi} e^{-(2N + \gamma + 2)}(2N + \gamma +
2)^{2N + \gamma + 3/2}\right\}}{\chi^{2N+\gamma+2}}
 \end{multline}
We can see now that the remainder is only \emph{algebraically} small
unless $N \sim \lvert \chi \rvert /2\epsilon$ (where the ratio of consecutive
terms are equal) and thus we set $N = r/2\epsilon + \rho$ where $\rho$ is
bounded as $\epsilon \to 0$. 

Although there are four homogeneous solutions to Equation (\ref{eq:remeq}) as
$\epsilon \to 0$, we will show that one in particular, 
\begin{equation}
 R_N(x) \sim S(x)e^{-i(x-\sigma)/\epsilon}
\end{equation}
is switched on as the Stokes line is crossed. We will call the
function $S(x)$ the \emph{Stokes multiplier}, and we expect it to vary
smoothly from one constant to another across the Stokes line. We write
\begin{equation} \label{eq:stokesvars}
\chi = x - \sigma = re^{i\theta} \text{\quad and \quad} \frac{d}{dx} =
-\frac{ie^{-i\theta}}{r} \frac{d}{d\theta},
\end{equation}
where now, since $N$ is fixed (and thus also the modulus, $r$), we
are only interested in the ``fast'' variation in $\theta$ across the Stokes
line. Then using Equations (\ref{eq:rhs}) and (\ref{eq:stokesvars}) in
(\ref{eq:remeq}) gives
\begin{multline} \label{eq:remainmessy}
\frac{dS}{d\theta} \sim
\frac{\Lambda \sqrt{r\pi}}{\sqrt{2}\epsilon^{\gamma+1/2}} \
e^{-r/\epsilon}e^{ire^{i\theta}/\epsilon} \Bigl( e^{-i\theta}
\Bigr)^{r/\epsilon + 2\rho + \gamma + 2}\Bigl( e^{-\pi i/2}\Bigr)^{r/\epsilon
+ 2\rho} e^{i\theta} \\
= \frac{\Lambda \sqrt{r\pi}}{\sqrt{2}\epsilon^{\gamma+1/2}} \times
\exp\left[-\frac{r}{\epsilon}\left\{1 - ie^{i\theta} + i\theta + \frac{\pi i}{2}
\right\} \right. \\ 
\left. + i\left\{ -2\rho\left(\theta+\frac{\pi}{2}\right) - \theta(\gamma+1)
\right\}\right]
\end{multline}
From the terms within the curly braces, we see that
the change in $S$ is exponentially small, except near the Stokes line $\theta =
-\pi/2$. Here, we will re-scale $\theta = -\pi/2 + \sqrt{\epsilon}\eta$ and
integrate Equation (\ref{eq:remainmessy}) from left ($\eta \to -\infty$) to
right to show that
\begin{equation}
S \sim \text{const} + \frac{\Lambda\sqrt{\pi}}{\sqrt{2}\epsilon^\gamma} e^{\pi
i(\gamma+1)/2}\int_{-\infty}^{\sqrt{r}\eta} e^{-s^2/2} \ ds.
\end{equation}
This integral (the error function) precisely illustrates
the \emph{smoothing} of the Stokes line in Figure \ref{fig}. Thus the jump in
the Stokes multiplier and consequently, the remainder is
\begin{equation} \label{eq:remainder}
\Bigl[ S \Bigr]_\text{Stokes}
\sim \frac{\Lambda \pi}{\epsilon^\gamma}e^{3\pi i/2} \Longrightarrow \Bigl[ R_N
\Bigr]_\text{Stokes} \sim
\frac{\Lambda \pi}{\epsilon^2} e^{3\pi i/2} e^{-i(x-\sigma)/\epsilon}.
\end{equation}
We must remember that the analysis must be repeated for analytic
continuation into the lower-half $x$-plane and thus near the singularity at $x =
-\pi i/2\gamma$. The result is another exponentially small contribution which is
the complex conjugate of Equation (\ref{eq:remainder}) and thus along the
real axis, the sum of contributions from crossing the pair of Stokes lines is, 
\begin{equation} \label{eq:fin_exp}
u_\text{exp} \sim -\frac{2\Lambda \pi}{\epsilon^2} e^{-\pi/2\gamma\epsilon}
\sin\left( x/\epsilon\right).
\end{equation}
Let us recap our analysis: (1) The \emph{singular} nature of the 5KdV
equation produces singularities in the early terms, (2) As more and more terms
are taken, the effects of the singularities grow, eventually producing
factorial over power divergence in the late terms, (3) Stokes lines emerge from
each of the singularities, and (4) By optimally truncation and examining the
jump in the remainder as the Stokes lines are crossed, we see the Stokes
Phenomenon and thus the appearance of exponentially small terms. 

So finally, we are ready to answer the original question: \emph{Do
there exist classical solitary wave solutions of the 5KdV equation?} No. For
suppose that we did impose the condition that only the base
(non-oscillatory) asymptotic solution applies at $x = -\infty$. Then as we pass
through $x = 0$, the term in Equation (\ref{eq:fin_exp}) necessarily switches on
and $u \sim u_0 + u_\text{exp}$ for $x > 0$. 

We have thus passed through Stokes' mist and subsequently, realized that
\emph{there do not exist classical solitary wave solutions of the 5KdV}.

\end{document}